\documentclass[10pt]{amsart}
\usepackage{amsmath,amssymb}
\newcommand{\doublespace}
   {\addtolength{\baselineskip}{0.15\baselineskip}}

\newtheorem{pdef}{Definition}[section] %
\newtheorem{thm}[pdef]{Theorem}        
\newtheorem{cor}[pdef]{Corollary}
\newtheorem{lem}[pdef]{Lemma}
\newtheorem{prop}[pdef]{Proposition}

\newcounter{equationnumber}
\renewcommand{\theequation}{\thesection.\arabic{equation}}
\def\mathletters{
    \addtocounter{equation}{1}
    \edef\@currentlabel{\theequation}
    \setcounter{equationnumber}{\value{equation}}
    \setcounter{equation}{0}
    \edef\theequation{\@currentlabel\noexpand\alph{equation}}
    }

\title{Closability property of operator algebras generated by normal operators
and operators of class $C_0$}
\author{Hao-Wei Huang}
\address{Department of Mathematics, Indiana University, 831 East 3rd Street,
Bloomington, IN 47405}
\email{huang39@indiana.edu}

\begin{document}
\maketitle \doublespace
\pagestyle{myheadings} \thispagestyle{plain}
\markboth{   }{ }

\begin{abstract}

An operator algebra $\mathcal{A}$ acting on a Hilbert space is said
to have the closability property if every densely defined linear
transformation commuting with $\mathcal{A}$ is closable. In this
paper we study the closability property of the von Neumann algebra
consisting of the multiplication operators on $L^2(\mu)$, and give
necessary and sufficient conditions for a normal operator $N$ such
that the von Neumann algebra generated by $N$ has the closability
property. We also give necessary and sufficient conditions for an
operator $T$ of class $C_0$ such that the algebra generated by $T$
in the weak operator topology and the algebra
$H^\infty(T)=\{u(T):u\in H^\infty\}$ have the closability property.
\end{abstract}
 \footnotetext[2]{{\it Key
words and phrases.}\, closability property, multiplication operator,
normal operator, operator of class $C_0$, Jordan model,
quasisimilarity.}

\section{Introduction}

An operator algebra $\mathcal{A}$ with the property that every
densely defined linear transformation in the commutant of
$\mathcal{A}$ is closable is said to have the closability property.
The closability property problem which has close connection with
transitive algebra problem was first studied by W. Arveson. In [4],
Arveson showed that the algebra $L^\infty$ acting on the Hilbert
space $L^2$ and the algebra $H^\infty$ acting on the Hardy space
$H^2$ have the closability property. In [2], H. Bercovici, R.G.
Douglas, C. Foias, and C. Pearcy showed that algebras
$\mathcal{W}_S$ and $\mathcal{W}_{S(\theta)}$ generated by the
unilateral shift $S$ and the Jordan block $S(\theta)$ (see Section 2
for precise definitions) in the weak operator topology,
respectively, and any maximal abelian selfadjoint subalgebra with a
cyclic vector have the closability property. They introduced some
general viewpoints, such as rationally strictly cyclic vector and
confluence for an algebra $\mathcal{A}$, to determine whether
$\mathcal{A}$ has the closability property. They also showed that if
an algebra $\mathcal{A}_1$ is a quasiaffinie transform of an algebra
$\mathcal{A}_2$ which has the closability property then
$\mathcal{A}_1$ has the closability property as well. As a
consequence, every unital commutative algebra
$\mathcal{A}\subset\mathcal{B(H)}$ with a rationally strictly cyclic
vector and the commutant of any operator of class $C_0$ have the
closability property. In particular, the algebra
$H^\infty(S(\theta)^*)=\{u(S(\theta)^*):u\in H^\infty\}$ has the
closability property for any nonconstant inner function $\theta$
which was proved independently by D. Sarason (see [7]). We refer the
reader to [2] for further details about confluent operator algebras
and the effect of quasimilarity on the closability property.

In [2, Proposition 3.5], some examples of operator algebras without
the closability property were given which point out that an algebra
with the closability property must be sufficiently large and should
not have uniform infinite multiplicity. This observation motivates
us investigate the relation between the closability property and
uniform finite multiplicity in this paper.

This paper is organized as follows. In Section 2 we give the
preliminaries and terminology of multiplication operators and
operators of class $C_0$. In Section 3 we study the closability
property of the von Neumann algebras $\mathcal{A}_\mu$ consisting of
the multiplication operators on $L^2(\mu)$ and $\mathcal{W}^*_N$
generated by the normal operator $N$. We show that
$\mathcal{A}_\mu^{(n)}$ has the closability property if and only $n$
is finite while $\mathcal{W}^*_N$ has the closability property if
and only if $N$ has uniform finite multiplicity. In the study of the
closability property, it is essential to examine closed unbounded
linear transformations in the commutant of a bounded operator. In
Section 4 we characterize the closed, densely defined linear
transformations intertwining two operators of class $C_0$. In
Section 5 we deal with the closability property of unital algebras
$H^\infty(T)=\{u(T):u\in H^\infty\}$ and $\mathcal{W}_T$ which is
generated by $T$ in the weak operator topology where $T$ is an
operator of class $C_0$, and show that $H^\infty(T)$ has the
closability property if $T$ has finite multiplicity. We also provide
necessary and sufficient conditions for $T$ with infinity
multiplicity such that the algebra $H^\infty(T)$ has the closability
property. Moreover, we show that the algebra $\mathcal{W}_T$ has the
closability property if and only if the algebra $H^\infty(T)$ does.

The author wishes to thank his advisor, Hari Bercovici, for his
generosity, his mathematical insight and for suggesting the problems
investigated in this paper.

\section{Preliminary}

Throughout this paper, the Hilbert space is over the complex number
$\mathbb{C}$. The space of bounded linear operators
$T:\mathcal{H}\rightarrow\mathcal{K}$, where $\mathcal{H}$ and
$\mathcal{K}$ are Hilbert spaces is denoted by $\mathcal{B(H,K)}$.
We will write $\mathcal{B(H)}$ instead of $\mathcal{B(H,H)}$, and
denote the range and kernel space of an operator $T$ by
$\mathrm{ran}\;T$ and $\mathrm{ker}\;T$, respectively. If
$\mathcal{M}$ is a submanifold of $\mathcal{H}$ then
$\overline{\mathcal{M}}$ is the norm closure of $\mathcal{M}$ and
$\mathcal{M}^\perp$ is the orthogonal complement of $\mathcal{M}$.
Denote by $P_\mathcal{M}$ the orthogonal projection of $\mathcal{H}$
onto $\mathcal{M}$ when $\mathcal{M}$ is closed. For an arbitrary
subalgebra $\mathcal{A}$ of $\mathcal{B(H)}$ and any collection
$\mathcal{L}$ of closed subspaces of $\mathcal{H}$,
$\mathrm{Lat}(\mathcal{A})$ means the lattice of invariant subspaces
of $\mathcal{A}$ while $\mathrm{Alg}(\mathcal{L})$ means the set of
those $A\in\mathcal{B(H)}$ such that
$A\mathcal{M}\subset\mathcal{M}$ for all
$\mathcal{M}\in\mathcal{L}$. For any algebra
$\mathcal{A}\subset\mathcal{B(H)}$ and $1\leq n\leq\infty$, define
the algebra $\mathcal{A}^{(n)}$ by
\[\mathcal{A}^{(n)}=\{T^{(n)}=\underbrace{T\oplus\cdots\oplus T}_{n\;\;\mathrm{times}}:T\in\mathcal{A}\}\]
which acts on the Hilbert space
\[\mathcal{H}^{(n)}=\underbrace{\mathcal{H}\oplus\cdots\oplus\mathcal{H}}_{n\;\;\mathrm{times}}.\] Given
an operator $T$ in $\mathcal{B(H)}$, the unital algebras generated
by $T$ in the weak operator topology and in the weak$^*$ topology
are denoted by $\mathcal{W}_T$ and $\mathcal{A}_T$, respectively.
For any subset $\mathcal{S}\subset\mathcal{B(H)}$, $\mathcal{S}'$ is
the set of operators commuting with every elements of $\mathcal{S}$
and called the commutant of $\mathcal{S}$. A von Neumann algebra
$\mathcal{A}$ is a unital $C^*$-algebra contained in
$\mathcal{B(H)}$ that is closed in the weak operator topology, i.e.,
a unital $C^*$-subalgebra of $\mathcal{B(H)}$ such that
$\mathcal{A}''=\mathcal{A}$ in light of the double commutant
theorem. The von Neumann algebra generated by $T$, i.e., the
smallest von Neumann algebra containing $T$ will be denoted by
$\mathcal{W}^*_T$. An important consequence of the double commutant
theorem is the equality
\[\mathcal{W}^*_N=\{N\}'',\]
where $N$ is a normal operator.

If $\mu$ is a compactly supported, regular Borel measure on
$\mathbb{C}$ and $N_\mu:L^2(\mu)\rightarrow L^2(\mu)$ is defined by
$(N_\mu f)(z)=zf(z)$, $f\in L^2(\mu)$, then $N_\mu$ is a normal
operator. Given a $\sigma$-finite measure space
$(\mathcal{X},\Omega,\mu)$ and function $\phi\in L^\infty(\mu)$, let
$M_\phi$ be the multiplication operator $M_\phi f=\phi f$ on
$L^2(\mu)$. If $\mathcal{A}_\mu=\{M_\phi:\phi\in L^\infty(\mu)\}$
then $\mathcal{A}_\mu$ is a von Neumann algebra with the property
$\mathcal{A}_\mu=\mathcal{A}_\mu'=\mathcal{A}_\mu''$. If $\mu$ is
assumed to be compactly supported on $\mathbb{C}$ then
$\{N_\mu\}'=\mathcal{A}_\mu$. If $\mathcal{H}$ is separable and
$N\in\mathcal{B(H)}$ is a normal operator then there exist mutually
singular measures $\mu_\infty,\mu_1,\mu_2,\cdots$ (some of which may
be zero measures) such that $N$ is unitarily equivalent to
\[N_{\mu_\infty}^{(\infty)}\oplus N_{\mu_1}\oplus N_{\mu_2}^{(2)}\oplus\cdots.\]
If $\Delta_n$ denotes the support of the measure $\mu_n$,
$n=1,2,\cdots,\infty$, the function
$m_N:\mathbb{C}\rightarrow\{0,1,\cdots,\infty\}$ associated with the
normal operator $N$ defined by
$m_N=\infty\cdot\chi_{\Delta_\infty}+\chi_{\Delta_1}+2\chi_{\Delta_2}+\cdots$
is a Borel function and called the multiplicity function of $N$. A
normal operator is said to have uniform finite multiplicity if its
multiplicity function is finite a.e. on $\mathbb{C}$.

Denote by $H^p$, $0<p\leq\infty$, the usual Hardy spaces on the unit
disc $\mathbb{D}$. For two functions $\theta$ and $\theta'$ in
$H^\infty$, we say that $\theta$ divides $\theta'$ or $\theta'$ is a
divisor of $\theta$, denoted by $\theta|\theta'$, if
$\theta'=\theta\phi$ for some $\phi\in H^\infty$. If $\theta$ and
$\theta'$ differ only by a constant scalar factor of absolute value
one, i.e., $\theta|\theta'$ and $\theta'|\theta$, then we use the
notation $\theta\equiv\theta'$. For a family $F$ of functions in
$H^\infty$, the notation $\bigwedge F$ (or $\bigwedge_if_i$ if
$F=\{f_i:i\in I\}$, or $f_1\wedge f_2$ if $F=\{f_1,f_2\}$) stands
for the greatest common inner divisor of $F$. The least common inner
multiple of $F$ is denoted by $\bigvee F$ (or $\bigvee_if_i$ if
$F=\{f_i:i\in I\}$, or $f_1\vee f_2$ if $F=\{f_1,f_2\}$).

For any inner function $\theta$, the space
$\mathcal{H}(\theta)=H^2\ominus\theta H^2$ is an invariant subspace
for $S^*$, the adjoint of the unilateral shift $S$ on the Hardy
space $H^2$. The operator
$S(\theta)\in\mathcal{B}(\mathcal{H}(\theta))$ defined by the
requirement that $S(\theta)^*=S^*|\mathcal{H}(\theta)$ is called a
Jordan block. A contraction $T\in\mathcal{B(H)}$ is called
completely nonunitary if it does not have any nontrivial unitary
direct summand. For a completely nonunitary contraction $T$, the
Sz.-Nagy--Foias functional calculus is an algebra homomorphism
$u\mapsto u(T)\in\mathcal{B(H)}$ of the algebra $H^\infty$ which
extends the usual polynomial calculus. For instant given any $u\in
H^\infty$, $u(S)$ is the analytic Toeplitz operator on $H^2$, i.e.,
the multiplication operator by $u$, and
$u(S(\theta))=P_{\mathcal{H}(\theta)}u(S)|\mathcal{H}(\theta)$.

A completely nonunitary contraction $T\in\mathcal{B(H)}$ is said to
be of class $C_0$ if $u(T)=0$ for some nonzero function $u\in
H^\infty$. If $T$ is of class $C_0$, the ideal $\{u\in
H^\infty:u(T)=0\}$ is of the form $m_TH^\infty$, where $m_T$ is an
inner function, uniquely determined up to a constant factor of
absolute value one, and called the minimal function of $T$. The
operator $S(\theta)$ is one of operators of class $C_0$ and its
minimal function is $\theta$. One of the important things about
$S(\theta)$ is
\[\{S(\theta)\}'=\mathcal{W}_{S(\theta)}=\{u(S(\theta)):u\in
H^\infty\}.\] Any invariant subspace
$\mathcal{M}\subset\mathcal{H}(\theta)$ for $S(\theta)$ is of the
form $\mathcal{M}=\phi H^2\ominus\theta H^2$ where $\phi$ divides
$\theta$.

An operator $Q\in\mathcal{B}(\mathcal{H},\mathcal{H}')$ is a
quasiaffinity if it is an injection with dense range, i.e., if it
has a (possibly unbounded) inverse defined on a dense domain in
$\mathcal{H}'$. An operator $T\in\mathcal{B(H)}$ is called a
quasiaffinine transform of $T'\in\mathcal{B}(\mathcal{H}')$ if there
exists a quasiaffinity $Q\in\mathcal{B}(\mathcal{H},\mathcal{H}')$
satisfying $T'Q=QT$ which is denoted by $T\prec T'$. If $T\prec T'$
then $T'^*\prec T^*$. The operators $T$ and $T'$ are called
quasisimilar, denoted by $T\sim T'$, if $T\prec T'$ and $T'\prec T$.
If $T_1\sim T_2$ and $T_2\sim T_3$ then $T_1\sim T_3$. For an
arbitrary operator $T\in\mathcal{B(H)}$ the cyclic multiplicity
$\mu_T$ is the smallest cardinal of a subset $\mathcal{M\subset H}$
with the property that $\bigvee_{n=0}^\infty
T^n\mathcal{M}=\mathcal{H}$, where the symbol $\bigvee$ is closed
linear span, while $T$ is said to have finite cyclic multiplicity if
$\mu_T<\infty$. Quasisimilarity plays an important role in
classification of operators of class $C_0$. Every operator $T$ of
class $C_0$ is quasisimilar to a unique Jordan operator, i.e., to an
operator of the form
\[\bigoplus_{i}S(\theta_i)\]
where the values of $i$ are ordinal numbers and the inner functions
$\theta_i$ are subject to the conditions that $\theta_i\equiv1$ for
some $i\geq0$, $\theta_{i_2}|\theta_{i_1}$ whenever $i_1\leq i_2$,
and $\theta_{i_1}\equiv\theta_{i_2}$ if
$\mathrm{card}(i_1)=\mathrm{card}(i_2)$. The properties of operators
of class $C_0$ stated in this section are known in the literature.
For more details about such operators, the reader may consult [3]
and [7].

For reader's convenience, in the following theorem we state the
facts about operators of class $C_0$ that will be frequently used in
the sequel. We refer to [3, Proposition 2.4.9] for (1), [3,
Corollary 3.1.7] for (2), [3, Theorem 3.1.16] for (3), [3, Theorem
3.5.1] for (4) and (5), [3, Corollary 3.5.10] for (6), [3, Theorem
4.1.2], and [3, Corollary 4.1.6] for (7).

\begin{thm} Suppose that $T\in\mathcal{B(H)}$ and
$T'\in\mathcal{B}(\mathcal{H}')$ are operators of class $C_0$ and
that $\theta,\theta'\in H^\infty$ are two inner functions.
\\ (1) An operator of the form $v(T)$, $v\in H^\infty$, is injective if and only if $v\wedge
m_T\equiv1$, in which case $v(T)$ is a quasiaffinity.
\\ (2) For any inner function $\theta$, the adjoint
$S(\theta)^*$ is unitarily equivalent to $S(\theta^\sim)$, where
$\theta^\sim\in H^\infty$ is the function defined by
$\theta^\sim(z)=\overline{\theta(\overline{z})}$, $z\in\mathbb{D}$.
\\ (3) Suppose
$A\in\mathcal{B}(\mathcal{H}(\theta),\mathcal{H}(\theta'))$. Then
$S(\theta')A=AS(\theta)$ if and only if
$A=P_{\mathcal{H}(\theta')}u(S)|\mathcal{H}(\theta)$ where $u\in
H^\infty$ such that $\theta'|u\theta$.
\\ (4) We have
$T\prec T'$ if and only if $T'\prec T$.
\\ (5) If $T$ has finite cyclic multiplicity, then
$T\sim\bigoplus_{j=0}^{n-1}S(\theta_j)$ where inner functions
$\theta_0,\cdots,\theta_{n-1}$ satisfy the condition
$\theta_{j+1}|\theta_j$ for all $j$ and $m_T=\theta_0$.
\\ (6) Suppose that $\mathcal{H}$ and $\mathcal{H}'$ are separable and $\bigoplus_{j=0}^\infty
S(\theta_j)$ is the Jordan model of $T$. If $m_{T'}|\theta_j$ for
all $j$ then $\bigoplus_{j=0}^\infty S(\theta_j)$ is also the Jordan
model of $T\oplus T'$.
\\ (7) We have
\[\{T\}''=\{T\}'\cap\mathrm{AlgLat}(T)=\mathcal{A}_T=\mathcal{W}_T.\]
Moreover, there exists a function $v\in H^\infty$ with the
properties that $v(T)$ is a quasiaffinity and
$\mathcal{W}_T=\{v(T)^{-1}u(T)\in\mathcal{B(H)}:u\in H^\infty\}$.
\end{thm}

\section{Closability property and uniform finite multiplicity}

\begin{pdef} Suppose that $\mathcal{A}$ is a subalgebra of $\mathcal{B(H)}$. A densely
defined linear transformation with domain $\mathcal{D}(X)$ is said
to commute with the algebra $\mathcal{A}$ if $\mathcal{D}(X)$ is
dense in $\mathcal{H}$ and invariant under $\mathcal{A}$, and
$XT=TX$ on $\mathcal{D}(X)$ for any $T\in\mathcal{A}$, i.e., if its
graph
\[\mathcal{G}(X)=\{h\oplus Xh:h\in\mathcal{D}(X)\}\] is invariant
for $\mathcal{A}^{(2)}$.
\end{pdef}

Let $\mathcal{H}$ and $\mathcal{K}$ be Hilbert spaces and let $X$
and $Y$ be densely defined linear transformations from $\mathcal{H}$
into $\mathcal{K}$. Then $Y$ is called an extension of $X$ if
$\mathcal{G}(X)\subset\mathcal{G}(Y)$ and in symbols this is denoted
by $X\subset Y$. The linear transformation $X$ is called closed if
its graph is closed in $\mathcal{H}\oplus\mathcal{K}$, while it is
called closable if the closure $\overline{\mathcal{G}(X)}$ of its
graph is a graph of some linear transformation which is denoted by
$\overline{X}$ and called the closure of $X$. In fact, $X$ is
closable if and only if for any sequence $\{h_n\oplus Xh_n\}$ in
$\mathcal{G}(X)$ which converges to $0\oplus k$ as
$n\rightarrow\infty$ it follows that $k=0$. The following lemma
gives basic facts about certain products of bounded operators and
densely defined linear transformations.

\begin{lem} Let $A$ be any operator and $X$ a densely defined linear
transformation.
\\ (1) If $A$ is injective and $AX$ is closable then $X$ is
closable. Particularly, if $AX\subset B$ for some operator $B$ then
$X$ is closable.
\\ (2) If $XA$ is also a densely defined linear transformation and $X$ is closable, then $XA$ is closable.
\\ (3) Let $\mathcal{H}$ and $\mathcal{K}$ be Hilbert spaces and
set $W=P_{\mathcal{H}}|\mathcal{H}\oplus\mathcal{K}$. If the domain
of $X$ is dense in $\mathcal{H}$ and $XW$ is closable, then $X$ is
closable.
\end{lem}

\begin{pf} Pick any sequence $\{h_n\oplus Xh_n\}$ in
$\mathcal{G}(X)$ which converges to $0\oplus k$ as
$n\rightarrow\infty$. If $A$ is injective and $AX$ is closable then
$h_n\oplus AXh_n\in\mathcal{G}(AX)$ and $h_n\oplus
AXh_n\rightarrow0\oplus(Ak)$ as $n\rightarrow\infty$ which shows
that $Ak=0$, and so $k=0$. Note that $AX$ is closable if $AX\subset
B$ for some operator $B$. Hence (1) follows. Next, pick any sequence
$\{h_n\oplus XAh_n\}$ in $\mathcal{G}(XA)$ which converges to
$0\oplus k$ as $n\rightarrow\infty$. Since the sequence
$\{Ah_n\oplus XAh_n\}$ is in $\mathcal{G}(X)$ which converges to
$0\oplus k$ as $n\rightarrow\infty$, it follows that $k=0$ if $X$ is
closable. Hence (2) holds. By examining the graph
$\mathcal{G}(XW)=\mathcal{G}(X)\oplus\mathcal{K}$, it is easy to see
that (3) holds.
\end{pf} \qed

\begin{pdef} An unital algebra $\mathcal{A}$ of $\mathcal{B(H)}$ is
said to have the closability property if every densely defined
linear transformation commuting with $\mathcal{A}$ is closable.
\end{pdef}

Suppose that $(\mathcal{X},\Omega,\mu)$ is a measure space. Recall
that any function in $L^2(\mu)$ is the quotient of two bounded
functions. If $f\in L^2(\mu)$, for instance, define
\[v_f=\left\{
\begin{array}{ll}
1/f                     &\hbox{if\;\;\;$|f|>1$} \\
1                     &\hbox{if\;\;\;$|f|\leq1$} \\
\end{array}\;\;\;\;\;\mathrm{and}\;\;\;\;\;u_f=fv_f.
\right.\leqno{(3.1)}\] Then $u_f$ and $v_f$ are in $L^\infty(\mu)$,
$v_f$ is nonzero $\mu$-a.e., and $f=u_f/v_f$.

To investigate the main topic of this section, we need some lemmas.

\begin{lem} Suppose that $(\mathcal{X},\Omega,\mu)$ is a $\sigma$-finite measure space,
and let $\mathcal{H}=\bigoplus_{j=1}^n\mathcal{H}_j$ where
$\mathcal{H}_j=L^2(\mu)$, $j=1,2,\cdots,n$, and $n<\infty$. If $X$
is a densely defined linear transformation with dense domain
$\mathcal{D}(X)$ commuting with the von Neumann algebra
$\mathcal{A}_\mu^{(n)}$ then the manifold
\[\mathcal{D}(X)\cap[\mathcal{H}_1\oplus\cdots\oplus\mathcal{H}_{j-1}\oplus\{0\}\oplus\mathcal{H}_{j+1}\oplus\cdots
\oplus\mathcal{H}_n]\] is dense in
$\mathcal{H}_1\oplus\cdots\oplus\mathcal{H}_{j-1}\oplus\{0\}\oplus\mathcal{H}_{j+1}\oplus\cdots
\oplus\mathcal{H}_n$ for $j=1,2,\cdots,n$.
\end{lem}

\begin{pf} It suffices to show that the manifold
\[\mathcal{D}=\mathcal{D}(X)\cap[\mathcal{H}_1\oplus\cdots\oplus\mathcal{H}_{n-1}\oplus\{0\}]\]
is dense in
$\mathcal{H}_1\oplus\cdots\oplus\mathcal{H}_{n-1}\oplus\{0\}$. By
taking linear combinations of elements in $\mathcal{D}$, it is
equivalent to showing the claim that for any $j=1,2,\cdots,n-1$,
functions of the form $0\oplus\cdots\oplus0\oplus
h_j\oplus0\oplus\cdots\oplus0$, $h_j\in\mathcal{H}_j$, can be
approximated by elements of $\mathcal{D}$. Consequently, it is
enough to show that the claim holds when $j=1$. To show this, first
assume $\mu(\mathcal{X})<\infty$ and pick two sequences
$f_m=\bigoplus_{j=1}^nf_{j,m}$ and $g_m=\bigoplus_{j=1}^ng_{j,m}$ in
$\mathcal{D}(X)$ so that $f_m\rightarrow1\oplus0\oplus\cdots\oplus0$
and $g_m\rightarrow0\oplus\cdots\oplus0\oplus1$ in $\mathcal{H}$.
For any $m$, let $v_{f_{n,m}},u_{f_{n,m}},v_{g_{n,m}}$, and
$u_{g_{n,m}}$ be the functions defined as in (3.1). Then simple
computations show that
\begin{align*}
k_m:&=\left(\bigoplus_{j=1}^nM_{v_{f_{n,m}}u_{g_{n,m}}}\right)f_m-
\left(\bigoplus_{j=1}^nM_{u_{f_{n,m}}v_{g_{n,m}}}\right)g_m \\
&=\left(\bigoplus_{j=1}^{n-1}(v_{f_{n,m}}u_{g_{n,m}}f_{j,m}-u_{f_{n,m}}v_{g_{n,m}}g_{j,m})\right)\oplus0,
\end{align*}
and hence the vector $k_m$ belongs to $\mathcal{D}$ for all $m$
since its last term in the summand is $0$. Using the properties that
$v_{f_{n,m}},u_{g_{n,m}}\rightarrow1$ and $u_{f_{n,m}}\rightarrow0$
in $L^2(\mu)$ as $m\rightarrow\infty$, together with the fact that
$\|v_f\|_\infty,\|u_f\|_\infty\leq1$ for any $f\in L^2(\mu)$ yields
that $k_m$ converges to $1\oplus0\oplus\cdots\oplus0$ in
$\mathcal{H}$. Note that for any $\phi\in L^\infty(\mu)$,
$(\bigoplus_{j=1}^nM_\phi)k_m$ is a sequence in $\mathcal{D}$ which
converges to $\phi\oplus0\oplus\cdots\oplus0$ in $\mathcal{H}$ as
$m\rightarrow\infty$. Since $L^\infty(\mu)$ is dense in $L^2(\mu)$,
it is clear that the claim holds for $j=1$.

We next consider the general case, i.e., $\mathcal{X}$ is countable
union of disjoint sets $\{\Delta_l\}_{l=1}^\infty$ for which
$\mu(\Delta_l)<\infty$ for all $l$. Set
$\mathcal{K}_l=\bigoplus_{j=1}^nL^2(\mu|\Delta_l)$, $l\geq1$, where
the measure $\mu|\Delta_l$ gives the same values as $\mu$ does on
$\Delta_l$ and $0$ off $\Delta_l$ and any function in
$L^2(\mu|\Delta_l)$ is viewed as a function defined on $\mathcal{X}$
which vanishes off $\Delta_l$. With this identification, it is easy
to see that the identity
\[\mathcal{D}_l:=
\left(\bigoplus_{j=1}^nM_{\chi_{\Delta_l}}\right)\mathcal{D}(X)=\mathcal{D}(X)\cap\mathcal{K}_l\]
holds and that $\mathcal{D}_l$ is a dense manifold of
$\mathcal{K}_l$ for all $l$. This yields that $X|\mathcal{K}_l$ is
also a densely defined linear transformation with dense domain
$\mathcal{D}_l$ commuting with $\mathcal{A}_{\mu|\Delta_l}^{(n)}$,
and consequently
\[\mathcal{D}_l\cap\left[\left(\bigoplus_{j=1}^{n-1}L^2(\mu|\Delta_l)\right)\oplus\{0\}\right]\]
is dense in $(\bigoplus_{j=1}^{n-1}L^2(\mu|\Delta_l))\oplus\{0\}$ by
the first part of the proof. For any $\epsilon>0$ and any
$f\in(\bigoplus_{j=1}^{n-1}\mathcal{H}_j)\oplus\{0\}$, let $N$ be a
sufficiently large integer such that
$\|f-P_{\oplus_{l=1}^N\mathcal{K}_l}f\|_\mathcal{H}<\epsilon$. For
any $l=1,2,\cdots,N$, pick an element
$x_l\in\mathcal{D}_l\cap[(\bigoplus_{j=1}^{n-1}L^2(\mu|\Delta_l))\oplus\{0\}]$
so that $\|x_l-P_{\mathcal{K}_l}f\|_\mathcal{H}<\epsilon/N$. Then it
is clear that $x=\sum_{l=1}^Nx_l\in\mathcal{D}$ and
\[\|x-f\|_\mathcal{H}\leq\sum_{l=1}^N\|x_l-P_{\mathcal{K}_l}f\|_\mathcal{H}+
\|f-P_{\oplus_{l=1}^N\mathcal{K}_l}f\|_\mathcal{H}<2\epsilon.\]
Since $\epsilon$ is arbitrary, the lemma follows.
\end{pf} \qed

\begin{lem} Under the same assumptions of Lemma 3.4,
the manifold
\[\mathcal{D}(X)\cap[\{0\}\oplus\cdots\oplus\{0\}\oplus\mathcal{H}_j\oplus\{0\}\oplus\cdots\oplus\{0\}]\]
is dense in
$\{0\}\oplus\cdots\oplus\{0\}\oplus\mathcal{H}_j\oplus\{0\}\oplus\cdots\oplus\{0\}$
for $j=1,2,\cdots,n$.
\end{lem}

\begin{pf} It suffices to show that the assertion holds when $j=1$.
Set
$\mathcal{K}=\mathcal{H}_1\oplus\cdots\oplus\mathcal{H}_{n-1}\oplus\{0\}$,
and consider the linear transformation
$Y=P_{\mathcal{K}}X|\mathcal{K}:\mathcal{D}(Y)\subset\mathcal{K}\rightarrow\mathcal{K}$
with domain denoted by $\mathcal{D}(Y)$. Then $Y$ is densely defined
since $\mathcal{D}(Y)=\mathcal{D}(X)\cap\mathcal{K}$ is dense in
$\mathcal{K}$ by Lemma 3.5. Moreover, $Y$ commutes with the algebra
$\mathcal{A}_\mu^{(n-1)}\oplus\{0\}$. Indeed, for any $\phi\in
L^\infty(\mu)$ and $h\in\mathcal{D}(Y)$, we have
$((\bigoplus_{j=1}^{n-1}M_\phi)\oplus0)h=(\bigoplus_{j=1}^nM_\phi)h\in\mathcal{D}(Y)$
and
\begin{align*}
Y\left(\left(\bigoplus_{j=1}^{n-1}M_\phi\right)\oplus0\right)h&=P_\mathcal{K}X\left(\bigoplus_{j=1}^nM_\phi\right)h \\
&=P_\mathcal{K}\left(\bigoplus_{j=1}^nM_\phi\right)Xh \\
&=\left(\left(\bigoplus_{j=1}^{n-1}M_\phi\right)\oplus0
\right)P_\mathcal{K}Xh=\left(\left(\bigoplus_{j=1}^{n-1}M_\phi\right)\oplus0
\right)Yh.
\end{align*} Applying Lemma 3.4 to $Y$ implies that
\begin{align*}
&\mathcal{D}(Y)\cap[\mathcal{H}_1\oplus\cdots\oplus\mathcal{H}_{n-2}\oplus\{0\}\oplus\{0\}] \\
=&\mathcal{D}(X)\cap[\mathcal{H}_1\oplus\cdots\oplus\mathcal{H}_{n-2}\oplus\{0\}\oplus\{0\}]
\end{align*} is dense in
$\mathcal{H}_1\oplus\cdots\oplus\mathcal{H}_{n-2}\oplus\{0\}\oplus\{0\}$.
Continuing the same method up to finite times implies the desired
result. The proof is complete.
\end{pf} \qed

\begin{lem} Under the same assumptions of Lemma of 3.4, if
$X=0$ on the sets
$\mathcal{D}(X)\cap[\{0\}\oplus\cdots\oplus\{0\}\oplus\mathcal{H}_j\oplus\{0\}\oplus\cdots\oplus\{0\}]$,
$j=1,2,\cdots,n$, then $X\equiv0$ on $\mathcal{D}(X)$.
\end{lem}

\begin{pf} We will finish the proof by showing the claim that $X=0$ on the set
\[\mathcal{D}_m:=\mathcal{D}(X)\cap[\mathcal{H}_1\oplus\cdots\oplus\mathcal{H}_{m}\oplus\{0\}\oplus\cdots
\oplus\{0\}]\] for any $m=1,2,\cdots,n$. Obviously, the claim is
true if $m=1$. Suppose that the claim holds for $m\leq j$, where
$j\leq n-1$. Applying Lemma 3.5 to find two nonzero elements
$\Phi=\phi_1\oplus\cdots\oplus\phi_{j+1}\oplus0\oplus\cdots\oplus0$
and $\Psi=0\oplus\cdots\oplus0\oplus\psi\oplus0\oplus\cdots\oplus0$
in $\mathcal{D}(X)$ where $\psi\in\mathcal{H}_{j+1}$. We may assume
that $\phi_1,\cdots,\phi_{j+1},\psi$ are all bounded by multiplying
them by a bounded function. Then $X\Psi=0$ and
\begin{align*}
\left(\bigoplus_{k=1}^nM_\psi\right)X\Phi&=\left(\bigoplus_{k=1}^nM_\psi\right)
X\Phi-\left(\bigoplus_{k=1}^nM_{\phi_{j+1}}\right)X\Psi \\
&=X\left[\left(\bigoplus_{k=1}^nM_\psi\right)\Phi-\left(\bigoplus_{k=1}^nM_{\phi_{j+1}}\right)\Psi\right] \\
&=X[(\phi_1\psi)\oplus\cdots\oplus(\phi_j\psi)\oplus0\oplus\cdots\oplus0] \\
&=0
\end{align*} because $X=0$ on $\mathcal{D}_j$.
For any $\mu$-measurable set $\Delta$ with finite measure, making
use of Lemma 3.5 and the preceding identity yields that
\[\left(\bigoplus_{k=1}^nM_{\chi_\Delta}\right)X\Phi=0,\leqno{(3.2)}\] whence
\[X\Phi=0.\leqno{(3.3)}\] Indeed, first pick a sequence of functions $\{f_m\}$ in
$\mathcal{H}_{j+1}$ so that $0\oplus\cdots\oplus0\oplus
f_m\oplus0\oplus\cdots\oplus0$ belongs to $\mathcal{D}(X)$ for all
$m$ and $f_m\rightarrow\chi_\Delta$ in $L^2(\mu)$ as
$m\rightarrow\infty$ by Lemma 3.5. Note that
\begin{align*}
&\left(\bigoplus_{k=1}^nM_{\chi_\Delta
v_{f_m}}\right)(0\oplus\cdots\oplus0\oplus f_m\oplus0\oplus\cdots\oplus0) \\
&=0\oplus\cdots\oplus0\oplus\chi_\Delta
u_{f_m}\oplus0\oplus\cdots\oplus0\in\mathcal{D}(X)
\end{align*} and $\chi_\Delta
u_{f_m}$ converges to $\chi_\Delta$ as $m\rightarrow\infty$, where
functions $v_{f_m}$ and $u_{f_m}$ are defined as in (3.1). Then
replacing $\psi$ by $\chi_\Delta u_{f_m}$ and letting
$m\rightarrow\infty$ yield (3.2). Finally, for any
$h=h_1\oplus\cdots\oplus
h_{j+1}\oplus0\oplus\cdots\oplus0\in\mathcal{D}(X)$, pick a nonzero
function $v\in L^\infty(\mu)$ with the property that
$vh_1,\cdots,vh_{j+1}\in L^\infty(\mu)$. Then it follows from (3.3)
that
\[\left(\bigoplus_{m=1}^nM_v\right)Xh
=X\left[\left(\bigoplus_{m=1}^nM_v\right)h\right]=0,\] and
consequently $Xh=0$ by the injectivity of $\bigoplus_{m=1}^nM_v$.
Hence $X=0$ on $\mathcal{D}_{j+1}$ and the claim holds by induction.
This completes the proof.
\end{pf} \qed

It follows form [4, Lemma 3.2] that if $X$ is a densely defined
linear transformation commuting with $\mathcal{A}_\mu$ where $\mu$
is a probability measure then there exists an everywhere defined
measurable function $k$ such that $Xf=kf$ for every $f$ in the
domain of $X$. As a result, there exists a function $v\in
L^\infty(\mu)$ which is nonzero almost everywhere such that $a=vk\in
L^\infty(\mu)$ and $M_vX\subset M_a$. Combining the proof of Lemma
3.4 with the above conclusions, it is easy to extend Arveson's
results to $\sigma$-finite case.

\begin{thm} Suppose that $(\mathcal{X},\Omega,\mu)$ is a $\sigma$-finite measure space.
Then the von Neumann algebra $\mathcal{A}_\mu^{(n)}$ has the
closability property if and only if $n$ is finite.
\end{thm}

\begin{pf} The closability property of the
algebra $\mathcal{A}_\mu^{(n)}$ yields that it has uniform finite
multiplicity, i.e., $n$ is finite by [2, Proposition 3.5(5)].
Conversely, suppose $n<\infty$ and let $X$ be a densely defined
linear transformation with domain $\mathcal{D}(X)$ commuting with
$\mathcal{A}_\mu^{(n)}$. Set
$\mathcal{H}=\bigoplus_{j=1}^n\mathcal{H}_j$, where
$\mathcal{H}_j=L^2(\mu)$ for all $j$, and
\[\mathcal{D}_j=\mathcal{D}(X)\cap[\{0\}\oplus\cdots\oplus\{0\}\oplus
\mathcal{H}_j\oplus\{0\}\oplus\cdots\oplus\{0\}].\] Then Lemma 3.5
implies that there exist densely defined linear transformations
$X_{ij}:\mathcal{D}_j\subset\mathcal{H}_j\rightarrow\mathcal{H}_i$,
$i,j=1,2,\cdots,n$, such that
\[X(0\oplus\cdots\oplus0\oplus f_j\oplus0\oplus\cdots\oplus0)=(X_{1j}f_j)\oplus\cdots\oplus(X_{nj}f_j)\] for
any $0\oplus\cdots\oplus0\oplus
f_j\oplus0\oplus\cdots\oplus0\in\mathcal{D}_j$. It is apparent that
$X_{ij}$ commutes with the algebra $\mathcal{A}_\mu$ for
$i,j=1\cdots,n$. Indeed, if $\phi\in L^\infty(\mu)$ then
\begin{align*}
(M_\phi X_{1j}f_j)\oplus\cdots\oplus(M_\phi
X_{nj}f_j)&=\left(\bigoplus_{j=1}^nM_\phi\right)X(0\oplus\cdots\oplus0\oplus
f_j\oplus0\oplus\cdots\oplus0) \\
&=X\left(\bigoplus_{j=1}^nM_\phi\right)(0\oplus\cdots\oplus0\oplus f_j\oplus0\oplus\cdots\oplus0) \\
&=X(0\oplus\cdots\oplus0\oplus(\phi f_j)\oplus0\oplus\cdots\oplus0) \\
&=(X_{1j}M_\phi f_j)\oplus\cdots\oplus(X_{nj}M_\phi f_j),
\end{align*}
which shows that $M_\phi X_{ij}\subset X_{ij}M_\phi$. Let $v\in
L^\infty(\mu)$ be a nonzero function so that $M_vX_{ij}\subset
M_{a_{ij}}$, where $a_{i,j}\in L^\infty(\mu)$ for
$i,j=1,2,\cdots,n$, and set
\[A=(M_{a_{ij}}),\;\;\;\;\;\;B=\bigoplus_{j=1}^nM_v,\;\;\;\mathrm{and}\;\;\;\;Z=BX-A.\] Since $A$ and
$B$ are in the algebra $M_{n\times
n}(\mathcal{A}_\mu)=\{\mathcal{A}_\mu^{(n)}\}'$, it follows that $Z$
is a densely defined linear transformation with domain
$\mathcal{D}(X)$ commuting with $\mathcal{A}_\mu^{(n)}$. Moreover,
the preceding discussions show that $Z=0$ on the sets
$\mathcal{D}_1,\cdots,\mathcal{D}_n$ from which we deduce that $Z=0$
on $\mathcal{D}(X)$ or, equivalently, $BX\subset A$, by Lemma 3.6.
By virtue of Lemma 3.2(1), $X$ is closable since $B$ is injective,
and consequently $\mathcal{A}_\mu^{(n)}$ has the closability
property. The proof is complete.
\end{pf} \qed

Recall that an algebra
$\mathcal{A}_1\subset\mathcal{B}(\mathcal{H}_1)$ is a quasiaffine
transform of an algebra
$\mathcal{A}_2\subset\mathcal{B}(\mathcal{H}_2)$ if there exists a
quasiaffinity $Q\in\mathcal{B}(\mathcal{H}_1,\mathcal{H}_2)$ so that
for any $T_2\in\mathcal{A}_2$, we have $T_2Q=QT_1$ for some
$T_1\in\mathcal{A}_1$. It follows from [2, Proposition 5.2] that if
$\mathcal{A}_1$ and $\mathcal{A}_2$ are unital algebras,
$\mathcal{A}_1$ is a quasiaffine transformation of $\mathcal{A}_2$,
and if $\mathcal{A}_2$ has the closability property, then so does
$\mathcal{A}_1$. Particularly, this says that the closability
property is invariant under unitary equivalence.

\begin{thm} Suppose that $\mathcal{H}$ is a separable Hilbert space and
$N$ is a normal operator on $\mathcal{H}$. Then the von Neumann
algebra $\mathcal{W}^*_N$ generated by $N$ has the closability
property if and only if $N$ has uniform finite multiplicity.
\end{thm}

\begin{pf} First note that there exist mutually singular measures $\mu_\infty,\mu_1,\mu_2,\mu_3,\cdots$
(some of which may be zero) such that $N$ is unitarilly equivalent
to
\[N_{\mu_\infty}^{(\infty)}\oplus N_{\mu_1}\oplus N_{\mu_2}^{(2)}\oplus N_{\mu_3}^{(3)}\oplus\cdots,\]
and therefore
\[\mathcal{W}^*_N\cong
\mathcal{A}_{\mu_\infty}^{(\infty)}\oplus
\mathcal{A}_{\mu_1}\oplus\mathcal{A}_{\mu_2}^{(2)}\oplus\mathcal{A}_{\mu_3}^{(3)}\oplus\cdots.\]
Since the closability property is invariant under unitary
equivalence, by virtue of [2, Lemma 3.8] and Theorem 3.7 it follows
that $\mathcal{W}^*_N$ has the closability property if and only if
$N$ has uniform finite multiplicity.
\end{pf} \qed

\section{Unbounded linear maps intertwining operators of class $C_0$}

Recall that an operator
$A\in\mathcal{B}(\mathcal{H}_1,\mathcal{H}_2)$ is said to intertwine
operators $T_2\in\mathcal{B}(\mathcal{H}_2)$ and
$T_1\in\mathcal{B}(\mathcal{H}_1)$ if $AT_1=T_2A$. Denote by
$\mathcal{I}(T_1,T_2)$ the set of all operators in
$\mathcal{B}(\mathcal{H}_1,\mathcal{H}_2)$ intertwining $T_2$ and
$T_1$. Interest in this section is primarily looking at densely
defined linear transformations that intertwine operators of class
$C_0$ as the definition is given below.

\begin{pdef} Suppose that $T_1\in\mathcal{B}(\mathcal{H}_1)$ and $T_2\in\mathcal{B}(\mathcal{H}_2)$.
A densely defined linear transformation $X$ with domain
$\mathcal{D}(X)$ is said to intertwine operators $T_2$ and $T_1$ if
$\mathcal{D}(X)$ is dense in $\mathcal{H}_1$ and invariant under
$T_1$, and $T_2X=XT_1$ on $\mathcal{D}(X)$.
\end{pdef}

It is apparent that a densely defined linear transformation $X$
intertwines operators $T_2$ and $T_1$ if and only if its graph
$\mathcal{G}(X)$ is an invariant manifold for $T_1\oplus T_2$. Let
$J:\mathcal{H}_1\oplus\mathcal{H}_2\rightarrow\mathcal{H}_2\oplus\mathcal{H}_1$
be the isomorphism defined by $J(h_1\oplus h_2)=(-h_2)\oplus h_1$.
If $X$ intertwines operators $T_2$ and $T_1$ then $J\mathcal{G}(X)$
is invariant for $T_2\oplus T_1$ from which we infer that
$\mathcal{G}(X^*)=[J\mathcal{G}(X)]^\perp$ is invariant for
$T_2^*\oplus T_1^*$ or, equivalently, $X^*$ intertwines $T_1^*$ and
$T_2^*$. In addition, if $X$ is closable then its closure
$\overline{X}$ also intertwines $T_2$ and $T_1$. If $X$ intertwines
$T$ and $T$ for any $T$ in an algebra $\mathcal{A}$, then $X$
commutes with $\mathcal{A}$.

We next consider some special types of densely defined linear
transformations. Let $T_1\in\mathcal{B}(\mathcal{H}_1)$ and
$T_2\in\mathcal{B}(\mathcal{H}_2)$ be any two operators. If
$A\in\mathcal{I}(T_1,T_2)$ and $B\in\{T_1\}'$ is a quasiaffinity
then $AB^{-1}$ is a densely defined linear transformation
intertwining $T_2$ and $T_1$ since the graph
$\mathcal{G}(AB^{-1})=\{Bh_1\oplus Ah_1:h_1\in\mathcal{H}_1\}$ is
invariant for $T_1\oplus T_2$. Further, it is apparent that
$\mathcal{G}(AB^{-1})=\mathrm{ran}\;M$, where
\[M=\left(
\begin{array}{cc}
B & 0 \\
A & 0 \\
\end{array}
\right)\in\{T_1\oplus T_2\}'.\] Hence we have
\begin{align*}
\mathcal{G}((AB^{-1})^*)&=[J\mathcal{G}(AB^{-1})]^\perp \\
&=[J\mathrm{ran}\;M]^\perp=J[\mathrm{ran}\;M]^\perp=J\mathrm{ker}\;M^* \\
&=\{h_2\oplus h_1\in\mathcal{H}_2\oplus\mathcal{H}_1:A^*h_2=B^*h_1\} \\
&=\mathcal{G}(B^{*-1}A^*),
\end{align*} which shows that $(AB^{-1})^*=B^{*-1}A^*$. Therefore, $AB^{-1}$ is closable if and only if
$B^{*-1}A^*$ is densely defined in which case
$\overline{AB^{-1}}=(B^{*-1}A^*)^*$. In general, we have
\begin{align*}
J^*\mathcal{G}(B^{*-1}A^*)&=\{h_1\oplus
h_2\in\mathcal{H}_1\oplus\mathcal{H}_2:B^*h_1+A^*h_2=0\} \\
&=\mathrm{ker}\;M^*,
\end{align*} or, equivalently,
$\mathcal{G}(B^{*-1}A^*)=J\mathrm{ker}\;M^*$ whence
$[J^*\mathcal{G}(B^{*-1}A^*)]^\perp=\overline{\mathrm{ran}\;M}=\overline{\mathcal{G}(AB^{-1})}$.

If $T$ is an operator of class $C_0$ and $X$ is a closed, densely
defined linear transformation commuting with $T$ then
$X=\overline{AB^{-1}}$ for some operator $A\in\{T\}'$ and
quasiaffinity $B\in\{T\}'$ which was prove by Bervocivi [1]. The
following result, though not stated explicitly, is essentially due
to him.

\begin{prop} Let $T_1$ and $T_2$
be two operators of class $C_0$. If $X$ is a closed, densely defined
linear transformation intertwining $T_2$ and $T_1$ then
$X=\overline{AB^{-1}}=(B^{*-1}A^*)^*$, where
$A\in\mathcal{I}(T_1,T_2)$ and $B\in\{T_1\}'$ is a quasiaffinity.
Consequently, $\mathcal{G}(X)=\overline{\mathrm{ran}\;M}$, where
\[M=\left(
\begin{array}{cc}
B & 0 \\
A & 0 \\
\end{array}
\right)\in\{T_1\oplus T_2\}'.\]
\end{prop}

\begin{pf} Denote by $\mathcal{D}(X)$ the domain of $X$. Let
$T=T_1\oplus T_2|\mathcal{G}(X)$ and define the quasiaffinity
$Q:\mathcal{G}(X)\rightarrow\mathcal{H}_1$ by $Q(h_1\oplus
Xh_1)=h_1$, $h_1\in\mathcal{D}(X)$. It is clear that $T_1Q=QT$ by
the hypothesis, and hence $T\prec T_1$. This implies that
$TQ'=Q'T_1$ for some quasiaffinity
$Q'\in\mathcal{B}(\mathcal{H}_1,\mathcal{G}(X))$ by Theorem 2.1(4).
Assume $Q'h_1=Bh_1\oplus Ah_1$, $h_1\in\mathcal{H}_1$. Then it is
easy to see that $B\in\{T_1\}'$ and $T_2A=AT_1$. If $Bh_1=0$ for
some $h_1\in\mathcal{H}_1$, then $Q'h_1=0$ since $\mathcal{G}(X)$ is
a graph, and hence $h_1=0$. The fact that $Q'\mathcal{H}_1$ is dense
in $\mathcal{G}(X)$ yields that
$\mathcal{D}(X)\subset\overline{B\mathcal{H}_1}$, and hence $B$ has
dense range. Finally, the rest desired results follow from the
equalities
\[\overline{\mathcal{G}(AB^{-1})}=\overline{\{Bh\oplus Ah:h\in\mathcal{H}_1\}}=\overline{Q'
\mathcal{H}_1}=\mathcal{G}(X)\] and the fact that $AB^{-1}$ is
closable. The proof is complete.
\end{pf} \qed

By considering the adjoint $X^*$ of $X$ which intertwines operators
of class $C_0$, the preceding proposition implies the next result.

\begin{prop} Let $T_1\in\mathcal{B}(\mathcal{H}_1)$ and $T_2\in\mathcal{B}(\mathcal{H}_2)$
be any two operators of class $C_0$. If $X$ is a closed, densely
defined linear transformation intertwining $T_2$ and $T_1$ then
$X=B^{-1}A=(A^*B^{*-1})^*$, where $A\in\mathcal{I}(T_1,T_2)$ and
$B\in\{T_2\}'$ is a quasiaffinity. Consequently,
$\mathcal{G}(X)=J^*\mathrm{ker}\;M$, where
\[M=\left(
\begin{array}{cc}
B & A \\
0 & 0 \\
\end{array}
\right)\in\{T_2\oplus T_1\}'\] and $J$ is the isomorphism defined by
$J(h_1\oplus h_2)=(-h_2)\oplus h_1$, $h_1\in\mathcal{H}_1$,
$h_2\in\mathcal{H}_2$.
\end{prop}

\begin{pf} Since $\mathcal{G}(X^*)$ is an invariant subspace for
the operator $T_2^*\oplus T_1^*$ of class $C_0$, applying the
conclusions in Proposition 4.2 to the closed, densely defined linear
transformation $X^*$ yields that $X^*=\overline{A^*B^{*-1}}$ for
some $A^*\in\mathcal{I}(T_2^*,T_1^*)$ and quasiaffinity
$B^*\in\{T_2^*\}'$. Since $A^*B^{*-1}$ is closable, it follows that
$X=X^{**}=\overline{A^*B^{*-1}}^*=(A^*B^{*-1})^*=B^{-1}A$. The
identity $\mathcal{G}(X)=J^*\mathrm{ker}\;M$ is easy to verify. The
proposition follows.
\end{pf} \qed

\begin{cor} If $T$ is an operator of class $C_0$ and $X$ is a closed
densely defined linear transformation commuting with $T$ or the
algebra $\mathcal{W}_T$ then $X=B^{-1}A$ where $A,B\in\{T\}'$ and
$B$ is a quasiaffinity.
\end{cor}

\begin{pf} The first assertion follows immediately from Proposition 4.3. To
finish the proof, it suffices to show that $B^{-1}A$ also commutes
with $\mathcal{W}_T$. Based the facts that $\mathcal{W}_T=\{T\}''$
which is a part of Theorem 2.1(7) and
$\mathcal{D}(B^{-1}A)=\{h\in\mathcal{H}:Ah\in\mathrm{ran}\;B\}$, it
is easy to see that $\mathcal{G}(B^{-1}A)$ is an invariant subspace
for $\mathcal{W}_T^{(2)}$. The corollary follows.
\end{pf} \qed

Let $T_1$ and $T_2$ be any operators and suppose that
$A\in\mathcal{I}(T_1,T_2)$ and $B\in\{T_2\}'$ is injective. In
general, $B^{-1}A$ is a closed linear transformation that
intertwines $T_2$ and $T_1$ but not necessarily densely defined. If
$T_1$ and $T_2$ are of class $C_0$ then $B^{-1}A$ is densely defined
if and only if $BC=AD$ for some operator $C\in\mathcal{I}(T_1,T_2)$
and quasiaffinity $D\in\{T_1\}'$ by Proposition 4.2. Note that in
general the condition $BC=AD$ only implies that
$\overline{CD^{-1}}\subset B^{-1}A$. To precede further, we need the
following lemma to verify when two closed, densely defined linear
transformation are equal.

\begin{lem} Suppose that $T_1$ and $T_2$ are any operators and that
$A,C$ are in $\mathcal{I}(T_1,T_2)$ and $B,D\in\{T_2\}'$ are
quasiaffinities. If $EB=FD$ for some quasiaffinities
$E,F\in\{T_2\}'$ then $B^{-1}A=D^{-1}C$ if and only if $EA=FC$.
Particularly, if $BD=DB$ then $B^{-1}A=D^{-1}C$ if and only if
$BC=DA$.
\end{lem}

\begin{pf} Note that $B^{-1}A=(EB)^{-1}(EA)$ and
$D^{-1}C=(FD)^{-1}(FC)=(EB)^{-1}(FC)$. Then it is apparent that
$B^{-1}A=D^{-1}C$ if and only if $EA=FC$. If $BD=DB$ then letting
$E=D$ and $F=B$ yields the desired result.
\end{pf} \qed

For operators of class $C_0$ with finite multiplicity, we have the
following proposition.

\begin{prop} Suppose that $T_1$ and $T_2$
are any two operators of class $C_0$ with $\mu_{T_2}<\infty$. If $X$
is a closed, densely defined linear transformation intertwining
$T_2$ and $T_1$ then $X=v(T_2)^{-1}A$, where
$A\in\mathcal{I}(T_1,T_2)$ and $v\in H^\infty$ so that $v(T_2)$ is a
quasiaffinity.
\end{prop}

\begin{pf} By Proposition 4.3, $X=B^{-1}A_0$ for some
$A_0\in\mathcal{I}(T_1,T_2)$ and quasiaffinity $B\in\{T_2\}'$. Since
$T_2$ has finite multiplicity, it follows from [1, Proposition 2]
that $BC=v(T_2)$ where $C\in\{T_2\}'$ and $v\in H^\infty$ so that
$v(T_2)$ and $C$ are quasiaffinities. Let $A=CA_0$. Then the fact
that $B$ commutes with $v(T_2)$, together with the identity
$BA=v(T_2)A_0$, yields that $B^{-1}A_0=v(T_2)^{-1}A$ by Lemma 4.5.
This finishes the proof.
\end{pf} \qed

\begin{prop} Suppose that $T_1$ and
$T_2$ are operators of class $C_0$ and that
$A\in\mathcal{I}(T_1,T_2)$ and $v\in H^\infty$ so that
$v\wedge(m_{T_1}\vee m_{T_2})\equiv1$. Then $X=Av(T_1)^{-1}$ is a
closable, densely defined linear transformation intertwining
$u(T_2)$ and $u(T_1)$ for any $u\in H^\infty$ and its closure is
$\overline{X}=v(T_2)^{-1}A$.
\end{prop}

\begin{pf} It is clear that $X$ intertwines $u(T_2)$ and $u(T_1)$ for any $u\in H^\infty$, and that $v(T_2)X\subset A$
whence $X$ is closable by Theorem 2.1(1) and Lemma 3.2. Since
$\overline{X}$ intertwines $T_2$ and $T_1$, from Proposition 4.3 we
infer that $\overline{X}=D^{-1}C$ for some
$C\in\mathcal{I}(T_1,T_2)$ and quasiaffinity $D\in\{T_2\}'$ which
shows that $DA=Cv(T_1)=v(T_2)C$. Hence $\overline{X}=v(T_2)^{-1}A$
by Lemma 4.5. This completes the proof.
\end{pf} \qed

If the additional assumption $\mu_{T_1}<\infty$ is put in
Proposition 4.6, more results can be obtained.

\begin{prop} Suppose that $T_1$ and
$T_2$ are operators of class $C_0$ with finite multiplicity. If $X$
is a closed, densely defined linear transformation intertwining
$T_2$ and $T_1$ then $X=\overline{Av(T_1)^{-1}}=v(T_2)^{-1}A$ for
some $A\in\mathcal{I}(T_1,T_2)$ and $v\in H^\infty$ so that $v(T_1)$
and $v(T_2)$ are quasiaffinities.
\end{prop}

\begin{pf} The second equality in the assertion holds by Proposition 4.7. Next, suppose
$T_1\in\mathcal{B}(\mathcal{H}_1)$,
$T_2\in\mathcal{B}(\mathcal{H}_2)$ and denote by $\mathcal{D}(X)$ by
the domain of $X$. Then the closed, densely defined linear
transformation
\[Y=\left(
\begin{array}{cc}
0 & 0 \\
X & 0 \\
\end{array}
\right):\mathcal{D}(X)\oplus\mathcal{H}_2\rightarrow\mathcal{H}_1\oplus\mathcal{H}_2\]
commutes with $T_1\oplus T_2$ since the graph
$\mathcal{G}(Y)=\{(h_1\oplus h_2)\oplus(0\oplus
Xh_1):h_1\in\mathcal{D}(X),h_2\in\mathcal{H}_2\}$ is an invariant
subspace for $(T_1\oplus T_2)\oplus(T_1\oplus T_2)$. By virtue of
Proposition 4.6 and Theorem 2.1(1), we conclude that $Y=v(T_1\oplus
T_2)^{-1}M$ where $v\in H^\infty$ so that $v\wedge(m_{T_1}\vee
m_{T_2})\equiv1$ and $M\in\{T_1\oplus T_2\}'$. This shows that
$X=v(T_2)^{-1}A$ for some operator $A\in\mathcal{I}(T_1,T_2)$, as
desired.
\end{pf} \qed

For nonconstant inner functions $\theta$, $\phi$ and functions
$u,v\in H^\infty$ so that $\phi|u\theta$ and
$v\wedge(\theta\vee\phi)\equiv1$, the operator
$A=P_{\mathcal{H}(\phi)}u(S)|\mathcal{H}(\theta)$ is in
$\mathcal{I}(S(\theta),S(\phi))$ and $v(S(\theta))$, $v(S(\phi))$
are quasiaffinities by Theorem 2.1(1) and (2). Thus, by Proposition
4.7 and 4.8 we have proved the following corollary.

\begin{cor} Let $\theta$ and $\phi$ be nonconstant inner functions.
\\ (1) For any $u,v\in H^\infty$ with $\phi|u\theta$ and
$v\wedge(\theta\vee\phi)\equiv1$, the linear transformation
$X=P_{\mathcal{H}(\phi)}u(S)v(S(\theta))^{-1}$ is a closable,
densely defined linear transformation intertwining $u(S(\phi))$ and
$u(S(\theta))$ for any $u\in H^\infty$ and its closure is
\[\overline{X}=v(S(\phi))^{-1}P_{\mathcal{H}(\phi)}u(S)|\mathcal{H}(\theta).\]
(2) Any closed, densely defined linear transformation $X$
intertwining $S(\phi)$ and $S(\theta)$ is of the form
\[X=\overline{P_{\mathcal{H}(\phi)}u(S)v(S(\theta))^{-1}}=v(S(\phi))^{-1}P_{\mathcal{H}
(\phi)}u(S)|\mathcal{H}(\theta),\] where $v$ and $u$ are in
$H^\infty$ so that $v\wedge(\theta\vee\phi)\equiv1$ and
$\phi|u\theta$.
\end{cor}

We finish this section with the following theorem which will be used
in Section 5.

\begin{thm} Let $\theta$ and $\phi$ be nonconstant inner functions.
If $X$ is a densely defined linear transformation intertwining
$u(S(\phi))$ and $u(S(\theta))$ for any $u\in H^\infty$ then $X$ is
closable with closure
\[\overline{X}=v(S(\phi))^{-1}u(S(\phi))P_{\mathcal{H}(\phi)}|\mathcal{H}(\theta),\]
where $u,v\in H^\infty$ so that $\phi|u\theta$ and
$v\wedge(\theta\vee\phi)\equiv1$.
\end{thm}

\begin{pf} Let $\theta'\equiv\theta\vee\phi$, $\phi'=\theta'/\phi$,
and
$\mathcal{H}(\theta)^\perp=\mathcal{H}(\theta')\ominus\mathcal{H}(\theta)$.
Consider the linear transformation $Y=\phi'(S(\theta'))XW$ where
$W=P_{\mathcal{H}(\theta)}|\mathcal{H}(\theta')$. First note that
the domain of $Y$ is $\mathcal{D}(X)\oplus\mathcal{H}(\theta)^\perp$
which is dense in $\mathcal{H}(\theta')$ where $\mathcal{D}(X)$ as
before denotes the domain of $X$. Then for any $h\in\mathcal{D}(X)$,
$g\in\mathcal{H}(\theta)^\perp$, and $u'\in H^\infty$ we have
\begin{align*}
u'(S(\theta'))Y(h\oplus g)&=u'(S(\theta'))\phi'(S(\theta'))Xh \\
&=\phi'(S(\theta'))u'(S(\phi))Xh \\
&=\phi'(S(\theta'))Xu'(S(\theta))h \\
&=\phi'(S(\theta'))XP_{\mathcal{H}(\theta)}u'(S(\theta'))(h+g) \\
&=\phi'(S(\theta'))Yu'(S(\theta'))(h\oplus g)
\end{align*} where we use the fact that
$u'(S(\theta'))\phi'(S(\theta'))|\mathcal{H}(\phi)=\phi'(S(\theta'))u'(S(\phi))$.
Thus $Y$ commutes with the algebra $\mathcal{W}_{S(\theta')}$ and
further, by virtue of [2, Proposition 3.7] we derive that $Y$ is
closable. Since the operator $\phi'(S(\theta'))|\mathcal{H}(\phi)$
is injective, it is apparent that $X$ is closable by Lemma 3.2(1)
and (3). The desired form of the closure $\overline{X}$ follows from
Corollary 4.9(2).
\end{pf} \qed

\section{Closability property of algebras $\mathcal{W}_T$ and $H^\infty(T)$}

In this section, we will give necessary and sufficient conditions
for an operator $T$ of class $C_0$ so that the algebras
$\mathcal{W}_T$ and $H^\infty(T)=\{u(T):u\in H^\infty\}$ have the
closability property. We first prove some lemmas.

\begin{lem} Suppose that $\theta_0,\cdots,\theta_n$ are nonconstant inner
functions and $T=\bigoplus_{j=0}^nS(\theta_j)$. If $X$ is a densely
defined linear transformation with dense domain $\mathcal{D}(X)$
commuting with the algebra $H^\infty(T)$ then the manifold
\[\mathcal{D}(X)\cap[\{0\}\oplus\cdots\oplus\{0\}\oplus
\mathcal{H}(\theta_j)\oplus\{0\}\oplus\cdots\oplus\{0\}]\] is dense
in $\{0\}\oplus\cdots\oplus\{0\}\oplus
\mathcal{H}(\theta_j)\oplus\{0\}\oplus\cdots\oplus\{0\}$ for any
$j=0,1,\cdots,n$.
\end{lem}

\begin{pf} It is enough to show that the
manifold
\[\mathcal{D}_0=
\mathcal{D}(X)\cap[\mathcal{H}(\theta_0)\oplus\{0\}\oplus\cdots\oplus\{0\}]\]
is dense in
$\mathcal{H}(\theta_0)\oplus\{0\}\oplus\cdots\oplus\{0\}$. Let
$\mathcal{H}=\bigoplus_{j=0}^n\mathcal{H}(\theta_j)$. First observe
that the set of functions in $\bigoplus_{j=0}^nH^2$ whose
projections onto $\mathcal{H}$ belong to $\mathcal{D}(X)$ is a dense
manifold of $\bigoplus_{j=0}^nH^2$. Using this observation to find
sequences
$\{f_m^{(0)}\}_{m=1}^\infty,\{f_m^{(1)}\}_{m=1}^\infty,\cdots,\{f_m^{(n)}\}_{m=1}^\infty$
in $\bigoplus_{j=0}^n H^2$ so that
$P_\mathcal{H}f_m^{(i)}\in\mathcal{D}(X)$ for all $i,m$ and
\begin{align*}
f_m^{(i)}&:=f_m^{(i,0)}\oplus f_m^{(i,1)}\oplus\cdots\oplus
f_m^{(i,n)} \\
&\rightarrow e_i
\end{align*} as $m\rightarrow\infty$ for any $i=0,1,\cdots,n$, where $e_i$ is the standard basis in
$\bigoplus_{j=0}^nH^2$
with the constant function $1$ in the $i$-th entry and $0$
elsewhere. Let $\Phi_m=(f_m^{(i,j)})$ be the matrix composed of the
entries $f_m^{(i,j)}$ and let $a_m^{(i)}$ be the product of $(-1)^i$
and the determinant of the matrix obtained by deleting the first
column and $i$-th row of $\Phi_m$. Then it follows from linear
algebra that
\[\sum_{i=0}^na_m^{(i)}f_m^{(i)}=(\det\Phi_m)\oplus0\oplus\cdots\oplus0
\rightarrow1\oplus0\oplus\cdots\oplus0\leqno{(5.1)}\] in
$H^{1/(n+1)}$ as $m\rightarrow\infty$ which implies the greatest
common inner divisor of the inner parts of the functions
$\{\det\Phi_m\}$ is $1$. Indeed, if $h=\bigwedge_m(\det\Phi_m)$,
then (5.1) implies that the limit function $1\in hH^{1/(n+1)}$, and
hence $h\equiv1$. Next pick a bounded outer function $v$ so that the
function $vf_m^{(i,j)}$ is in $H^\infty$ for all $i,j$, and $m$.
Note that $v^na_m^{(i)}\in H^\infty$, and therefore we have that
$(v^na_m^{(i)})(T)P_\mathcal{H}f_m^{(i)}\in\mathcal{D}(X)$ for all
$i$ and $m$ by the hypothesis. Moreover, from (5.1) we deduce that
\begin{align*}
h_m:&=\sum_{i=0}^n(v^na_m^{(i)})(T)P_\mathcal{H}f_m^{(i)} \\
&=\sum_{i=0}^nP_\mathcal{H}v^n(T)(a_m^{(i)}f_m^{(i)}) \\
&=(P_{\mathcal{H}(\theta_0)}(v^n\det\Phi_m))\oplus0\oplus\cdots\oplus0\in\mathcal{D}_0
\end{align*}
for all $m$. Since
\[\bigwedge_{m=1}^\infty v^n\det\Phi_m\equiv1\leqno{(5.2)}\] by the preceding discussion, it follows that
\[P_{\mathcal{H}(\theta_0)}(v^n\det\Phi_m)\neq0\]
for some $m$, and consequently $h_m\neq0$ for some $m$. Hence
$\mathcal{D}_0$ contains a nonzero element from which we derive that
the closure
$\overline{\mathcal{D}_0}=\mathcal{H}(\theta_0)\oplus\{0\}\oplus\cdots\oplus\{0\}$.
Indeed, since the closure $\overline{\mathcal{D}_0}$ is an invariant
subspace for $S(\theta_0)\oplus0\oplus\cdots\oplus0$ there exists
some inner divisor $\phi_0$ of $\theta_0$ such that
\[\overline{\mathcal{D}_0}=(\phi_0H^2\ominus\theta_0H^2)\oplus\{0\}\oplus\cdots\oplus\{0\}\]
which shows that $\phi_0\equiv1$ by (5.2). This completes the proof.
\end{pf} \qed

\begin{lem} Under the same assumptions of Lemma 5.1, if $X=0$ on the sets
$\mathcal{D}(X)\cap[\{0\}\oplus\cdots\oplus\{0\}\oplus\mathcal{H}(\theta_j)\oplus\{0\}\oplus\cdots\oplus
\{0\}]$, $j=0,1,\cdots,n$, then $X\equiv0$ on $\mathcal{D}(X)$.
\end{lem}

\begin{pf} We will finish the proof by showing the claim that $X=0$ on
the sets
\[\mathcal{M}_m:=\mathcal{D}(X)\cap[\mathcal{H}(\theta_0)\oplus\mathcal{H}(\theta_1)\oplus\cdots
\oplus\mathcal{H}(\theta_m)\oplus\{0\}\oplus\cdots\oplus\{0\}],\;\;\;m=0,1,\cdots,n.\]
The case $m=0$ follows from the hypothesis. Suppose now that the
claim holds for $m\leq j$, where $j\leq n-1$. Using lemma 5.1 to
find two nonzero elements $f=P_\mathcal{H}(u_0\oplus\cdots\oplus
u_{j+1}\oplus0\oplus\cdots\oplus0)$ and
$g=P_\mathcal{H}(0\oplus\cdots\oplus0\oplus
u\oplus0\oplus\cdots\oplus0)$ in $\mathcal{D}(X)$, where
$u,u_0,\cdots,u_{j+1}\in H^2$. By multiplying functions
$u,u_0,\cdots,u_{j+1}$ by a bounded function in $H^2$, we may assume
$u,u_0,\cdots,u_{j+1}\in H^\infty$. Since $Xg=0$, it follows that
\begin{align*}
u(T)Xf&=u(T)Xf-u_{j+1}(T)Xg \\
&=X(u(T)f-u_{j+1}(T)g) \\
&= XP_\mathcal{H}[(uu_0)\oplus\cdots\oplus
(uu_j)\oplus0\oplus\cdots\oplus0] \\
&=0,
\end{align*} which yields that
\[XP_\mathcal{H}(u_0\oplus\cdots\oplus
u_{j+1}\oplus0\oplus\cdots\oplus0)=0\leqno{(5.3)}\] for any vector
$P_\mathcal{H}(u_0\oplus\cdots\oplus
u_{j+1}\oplus0\oplus\cdots\oplus0)\in\mathcal{D}(X)$ with
$u_0,\cdots,u_{j+1}\in H^\infty$ by Lemma 5.1. Finally, for any
vector $f=P_\mathcal{H}(u_0\oplus\cdots\oplus
u_{j+1}\oplus0\oplus\cdots\oplus0)\in\mathcal{D}(X)$ let $v$ be a
bounded outer function so that $vu_0,\cdots,vu_{j+1}$ are bounded.
Then (5.3) shows that
\begin{align*}
v(T)Xf&=Xv(T)f \\
&=XP_\mathcal{H}[(vu_0)\oplus\cdots\oplus(vu_{j+1})\oplus0\cdots\oplus0] \\
&=0.
\end{align*} Since $v$ is an outer function, $v(T)$ is a quasiaffinity,
and hence $Xf=0$. Then the claim is proved by induction and this
completes the proof.
\end{pf} \qed

We now can prove the first main result of this section.

\begin{thm} Let $\theta_0,\cdots,\theta_n$ be nonconstant
inner functions and let $T=\bigoplus_{j=0}^nS(\theta_j)$. Then the
algebra $H^\infty(T)$ has the closability property.
\end{thm}

\begin{pf} Let $\mathcal{H}=\bigoplus_{j=0}^n
\mathcal{H}(\theta_j)$ and $X$ be a densely defined linear
transformation with domain $\mathcal{D}(X)$ commuting with the
algebra $H^\infty(T)$. By virtue of Lemma 5.1, the set
\[\mathcal{D}_j:=
\mathcal{D}(X)\cap[\{0\}\oplus\cdots\oplus\{0\}\oplus\mathcal{H}(\theta_j)\oplus\{0\}\oplus\cdots\oplus\{0\}]\]
is dense in
$\{0\}\oplus\cdots\oplus\{0\}\oplus\mathcal{H}(\theta_j)\oplus\{0\}\oplus\cdots\oplus\{0\}$
for any $j=0,1,\cdots,n$. If each $\mathcal{D}_j$ is viewed as a
densely manifold of $\mathcal{H}(\theta_j)$, there exist densely
defined linear transformations
$X_{ij}:\mathcal{D}_j\subset\mathcal{H}(\theta_j)\rightarrow\mathcal{H}(\theta_i)$,
$0\leq i,j\leq n$, such that
\[X(0\oplus\cdots\oplus0\oplus
f_j\oplus0\oplus\cdots\oplus0)=\bigoplus_{i=0}^n(X_{ij}f_j)\] for
all $0\oplus\cdots\oplus0\oplus
f_j\oplus0\oplus\cdots\oplus0\in\mathcal{D}_j$. It is apparent that
$X_{ij}$ intertwines $u(S(\theta_i))$ and $u(S(\theta_j))$ for all
$u\in H^\infty$ and $i,j=0,\cdots,n$. Indeed, we have
$u(S(\theta_j))\mathcal{D}_j=u(T)\mathcal{D}_j\subset\mathcal{D}_j$
and
\begin{align*}
\bigoplus_{i=0}^n(X_{ij}u(S(\theta_j))f_j)&=X(0\oplus\cdots\oplus0\oplus
(u(S(\theta_j))f_j)\oplus0\oplus\cdots\oplus0) \\
&=Xu(T)(0\oplus\cdots\oplus0\oplus f_j\oplus0\oplus\cdots\oplus0) \\
&=u(T)X(0\oplus\cdots\oplus0\oplus
f_j\oplus0\oplus\cdots\oplus0) \\
&=\bigoplus_{i=0}^n(u(S(\theta_i))X_{ij}f_j),
\end{align*}
which implies that
\[X_{ij}u(S(\theta_j))f_j=u(S(\theta_i))X_{ij}f_j\] for any
$i,j=0,1\cdots,n$. In light of Theorem 4.10, there exist functions
$v_0,v_1,\cdots,v_n\in H^\infty$ and
$A_{ij}\in\mathcal{I}(S(\theta_j),S(\theta_i))$, $i,j=0,\cdots,n$,
with the properties that $v_i\wedge\theta_i\equiv1$ and
$v_i(S(\theta_i))X_{ij}\subset A_{ij}$ for any $i,j=0,1,\cdots,n$.
Let $A=(A_{ij})$, $B=\bigoplus_{j=0}^nv_j(S(\theta_j))$, and
$Z=BX-A$. Notice that $A,B\in\{T\}'$ and $B$ is injective. Then the
above arguments show that $Z$ is a densely defined linear
transformation with domain $\mathcal{D}(X)$ commuting with
$H^\infty(T)$ and vanishes on $\mathcal{D}_j$ for any $j$. It
follows from Lemma 5.2 that $Z=0$ on $\mathcal{D}(X)$, and hence
$BX\subset A$. Thus $X$ is closable by Lemma 3.2(1), and
consequently $H^\infty(T)$ has the closability property. This
completes the proof.
\end{pf} \qed

By Theorem 5.3, together with the fact that the closability property
is preserved under quasisimilarity, we have the following result.

\begin{cor} For any operator $T$ of class $C_0$ with finite
multiplicity, the algebra $H^\infty(T)$ has the closability
property.
\end{cor}

\begin{pf} Since $T$ is quasisimilar to its Jordan model
$\bigoplus_{j=0}^nS(\theta_j)$ where $n$ is finite and the algebra
$H^\infty(\bigoplus_{j=0}^nS(\theta_j))$ has the closability
property by Theorem 5.3, the algebra $H^\infty(T)$ has the
closability property as well by [2, Proposition 5.2(2)].
\end{pf} \qed

The following proposition is a direct consequence of Theorem 5.3 and
[2, Proposition 3.5(5)].

\begin{prop} For any nonconstant inner function $\theta$,
the algebra $H^\infty(S(\theta)^{(n)})$ has the closability property
if and only if $n<\infty$.
\end{prop}

We will make use of the following lemma to investigate the
closability property of $H^\infty(T)$ when $T$ has infinity
multiplicity.

\begin{lem} Suppose that $T_1$ and $T_2$ are completely nonunitary contraction. If the algebra $H^\infty(T_1\oplus T_2)$
has the closability property then so do the algebras $H^\infty(T_1)$
and $H^\infty(T_2)$.
\end{lem}

\begin{pf} Suppose $T_1\in\mathcal{B}(\mathcal{H}_1)$ and $T_2\in\mathcal{B}(\mathcal{H}_2)$ and let $\mathcal{H}
=\mathcal{H}_1\oplus\mathcal{H}_2$. For any densely defined linear
transformation $X$ with domain $\mathcal{D}(X)$ commuting with the
algebra $H^\infty(T_1)$,
\[Y=\left(
\begin{array}{cc}
X & 0 \\
0 & 0 \\
\end{array}\right):\mathcal{D}(X)\oplus\mathcal{H}_2\subset\mathcal{H}\rightarrow\mathcal{H}\]
is a densely defined linear transformation commuting with the
algebra $H^\infty(T_1\oplus T_2)$. If $H^\infty(T_1\oplus T_2)$ has
the closability property then $Y$ is closable whence $X$ is
closable. Indeed, if $X$ is regarded as a linear map from
$\mathcal{H}_1$ into $\mathcal{H}$ then
$Y=XP_{\mathcal{H}_1}|\mathcal{H}$, and so $X$ is closable by Lemma
3.2(3). Hence $H^\infty(T_1)$, as well as $H^\infty(T_2)$, has the
closability property.
\end{pf} \qed

\begin{thm} Suppose that $T$ is an operator of class $C_0$ acting on a separable Hilbert space.
If $\bigoplus_{j=0}^\infty S(\theta_j)$ is the Jordan model of $T$,
where $\theta_j$ is nonconstant for all $j$, then the algebra
$H^\infty(T)$ has the closability property if and only if
$\bigwedge_{j=0}^\infty\theta_j\equiv1$.
\end{thm}

\begin{pf} Let $\theta=\bigwedge_{j=0}^\infty\theta_j$. First assume that the algebra $H^\infty(T)$ has
the closability property. Since $\theta|\theta_j$ for any $j$, it
follows from Theorem 2.1(6) that $\bigoplus_{j=0}^\infty
S(\theta_j)$ is also the Jordan model of $T\oplus T'$ where
$T'=S(\theta)^{(\infty)}$ which shows that $T\sim T\oplus T'$. By
[2, Proposition 5.2(2)], the algebra $H^\infty(T\oplus T')$ also has
the closability property. But Lemma 5.6 yields that the algebra
$H^\infty(T')$ has the closability property, which is a
contradiction unless $\theta\equiv1$ by Proposition 5.5.

Conversely, by [2, Proposition 5.2(2)] we may assume that
$T=\bigoplus_{j=0}^\infty S(\theta)$ with $\theta\equiv1$ and let
$X$ be any densely defined linear transformation with domain
$\mathcal{D}(X)$ commuting with the algebra $H^\infty(T)$. For any
$m\geq0$, let $T_m=\bigoplus_{j=0}^mS(\theta_j)$ and
$\mathcal{H}_m=\bigoplus_{j=0}^m\mathcal{H}(\theta_j)$, and consider
the linear transformation
$Y_m=\theta_m(T)X\theta_m(T):\mathcal{D}(Y_m)\subset\mathcal{H}_m\rightarrow\mathcal{H}_m$
with domain
\[\mathcal{D}(Y_m)=\{h_m\in\mathcal{H}_m:\theta_m(T)h_m\in\mathcal{D}(X)\}.\]
Observe that if $h\in\mathcal{D}(X)$ then
$\theta_m(T)P_{\mathcal{H}_m}h=\theta_m(T)h\in\mathcal{D}(X)$, which
shows that $P_{\mathcal{H}_m}\mathcal{D}(X)\subset\mathcal{D}(Y_m)$
is dense in $\mathcal{H}_m$, and hence $Y_m$ is densely defined.
Moreover, $Y_m$ commutes with the algebra $H^\infty(T_m)$, and
consequently is closable by Theorem 5.3. Indeed, for any
$h_m\in\mathcal{D}(Y_m)$ and $u\in H^\infty$ we have
\[\theta_m(T)u(T_m)h_m=\theta_m(T)u(T)h_m=u(T)\theta_m(T)h_m\in\mathcal{D}(X)\]
which shows that $u(T_m)h_m\in\mathcal{D}(Y_m)$, and
\begin{align*}
Y_mu(T_m)h_m&=\theta_m(T)X\theta_m(T)u(T_m)h_m \\
&=\theta_m(T)Xu(T)\theta_m(T)h_m \\
&=u(T)\theta_m(T)X\theta_m(T)h_m \\
&=u(T_m)\theta_m(T)X\theta_m(T)h_m=u(T_m)Y_mh_m.
\end{align*}
Finally, suppose $\{f_n\oplus Xf_n\}$ is a sequence in
$\mathcal{G}(X)$ so that $f_n\oplus Xf_n\rightarrow 0\oplus g$ as
$n\rightarrow\infty$. Then for any fixed $m$ we have
$P_{\mathcal{H}_m}f_n\rightarrow0$ and
\[Y_mP_{\mathcal{H}_m}f_n=\theta_m(T)X\theta_m(T)f_n=\theta_m(T)^2Xf_n\rightarrow\theta_m(T)^2g\] as
$n\rightarrow\infty$, which implies that $\theta_m(T)^2g=0$ for any
$m\geq0$ or, equivalently,
\[\theta_j|\theta_m^2g_j\;\;\;\;\mathrm{for}\;\;\;\mathrm{any}\;\;\;m,j
\geq0\leqno{(5.4)}\] if $g=\bigoplus_{j=0}^\infty g_j$,
$g_j\in\mathcal{H}(\theta_j)$. It follows that $\theta_j|g_j$ for
all $j$, i.e., $g=0$. Indeed, for any fixed $n,j\geq0$ and for any
$m\geq\max\{n,j\}$, by (5.4) we have
$(\theta_j/\theta_m)|\theta_ng_j$, and so $\bigvee_{m\geq
n,j}(\theta_j/\theta_m)|\theta_ng_j$ which is equivalent to
$\theta_j|\theta_ng_j$ since $\bigvee_{m\geq
n,j}(\theta_j/\theta_m)\equiv\theta_j/(\bigwedge_{m\geq
n,j}\theta_m)\equiv\theta_j$. Hence we proved that
$\theta_j|\theta_ng_j$ for any $n,j\geq0$. Repeating the above
argument yields the desired result that $\theta_j|g_j$ for any $j$.
Therefore, $X$ is closable, and so the algebra $H^\infty(T)$ has the
closability property. This completes the proof.
\end{pf} \qed

Assume now that $T\in\mathcal{B(H)}$ is an operator of class $C_0$
with the Jordan model $\bigoplus_iS(\theta_i)$. If
\[\bigwedge_{j<\omega}\theta_j\equiv1\leqno{(5.5)}\] then
$\theta_\omega\equiv1$, and so $\bigoplus_i
S(\theta_i)=\bigoplus_{j<\omega}S(\theta_j)$ and the underlying
Hilbert space $\mathcal{H}$ is separable. Therefore, it follows from
Theorem 5.7 that $H^\infty(T)$ has the closability property if (5.5)
holds. Conversely, if $\theta=\bigwedge_{j<\omega}\theta_j$ then
Theorem 2.1(7) shows that
\[\left(\bigoplus_{j<\omega}S(\theta_j)\right)\oplus
S(\theta)^{(\infty)}\sim\bigoplus_{j<\omega}S(\theta_j),\] and hence
by comparison of Jordan models we have $T\oplus
S(\theta)^{(\infty)}\sim T$. Hence it is easy to see from [2,
Proposition 5.2(2)], Proposition 5.5, and Lemma 5.6 that (5.5) holds
if the algebra $H^\infty(T)$ has the closability property. We
summarize these conclusions in the following theorem.

\begin{thm} Suppose that $T\in\mathcal{B(H)}$ is an operator of class $C_0$ with the Jordan model
$\bigoplus_iS(\theta_i)$. Then the algebra $H^\infty(T)$ has the
closability property if and only if
$\bigwedge_{j<\omega}S(\theta_j)\equiv1$ in which case $\mathcal{H}$
is separable.
\end{thm}

\begin{cor} For any  operator $T$ of
class $C_0$, the algebra $H^\infty(T)$ has the closability property
if and only if the algebra $H^\infty(T^*)$ has the closability
property.
\end{cor}

\begin{pf} If $H^\infty(T)$ has the closability property, the Jordan model of $T$ must be of the form
$\bigoplus_{j=0}^\infty S(\theta_j)$ with
$\bigwedge_{j=0}^\infty\theta_j\equiv1$ by Theorem 5.8. Note that we
have $T^*\sim\bigoplus_{j=0}^\infty S(\theta_j^\sim)$ by Theorem
2.1(2) and $\bigwedge_j\theta_j^\sim\equiv1$. As a consequence of
Theorem 5.8 and [2, Proposition 5.2(2)], the algebra $H^\infty(T^*)$
has the closability property. The same is also true for the converse
and the corollary is proved.
\end{pf} \qed

Recall that an operator $T$ of class $C_0$ is said to have the
property $(P)$ if every injection $A\in\{T\}'$ is a quasiaffinity.
It was shown in [3, Theorem 7.1.9] that $T$ has the property $(P)$
if and only if its Jordan model satisfies the condition (5.5). We
next investigate the closability property of other algebras
generated by an operator of class $C_0$. Recall from Theorem 2.1(7)
that
\[\{T\}''=\{T\}'\cap\mathrm{AlgLat}(T)=\mathcal{A}_T=\mathcal{W}_T=\mathcal{F}_T\] for any operator $T$
of class $C_0$. By virtue of [2, Lemma 3.3], if $H^\infty(T)$ has
the closability property then the algebra $\mathcal{W}_T$ as well as
$\mathrm{AlgLat}(T)$ has the closability property (it is shown in
[2] that the commutant $\{T\}'$ of any operator $T$ of class $C_0$
has closability property). The following theorem states that the
converse is also true.

\begin{thm} For any operator $T$ of class $C_0$ with Jordan model $\bigoplus_iS(\theta_i)$,
the followings are equivalent:
\\ (1) the algebra $H^\infty(T)$ has the closability property;
\\ (2) $\bigwedge_{j<\omega}\theta_j\equiv1$;
\\ (3) $T$ has the property $(P)$;
\\ (4) the algebra $\mathcal{W}_T$ has the closability
property.
\end{thm}

\begin{pf} The equivalences of (1), (2), and (3) were proved. By [2, Lemma 3.3], it suffices to show that (4)
implies (1), i.e., any densely defined linear transformation $X$
commuting with the algebra $H^\infty(T)$ is closable if the algebra
$\mathcal{W}_T$ has the closability property. Suppose
$T\in\mathcal{B(H)}$ and define
$\mathcal{D}=\{Ah:A\in\mathcal{W}_T,\;h\in\mathcal{D}(X)\}$ where
$\mathcal{D}(X)$ is the domain of $X$. By Theorem 2.1(7), let $v$ be
in $H^\infty$ with the properties that $v(T)$ is a quasiaffinity and
$\mathcal{W}_T=\{v(T)^{-1}u(T)\in\mathcal{B(H)}:u\in H^\infty\}$.
Then $\mathcal{D}(X)\subset\mathcal{D}$ is dense in $\mathcal{H}$
and the domain of the linear transformation $Xv(T)$ contains
$\mathcal{D}$. Indeed, for any $Ah\in\mathcal{D}$ we have
$v(T)(Ah)=(Av(T))h\in\mathcal{D}(X)$ since $Av(T)\in H^\infty(T)$
which shows that $Ah\in\{k\in\mathcal{H}:v(T)k\in\mathcal{D}(X)\}$,
the domain of $Xv(T)$. Moreover, the densely defined linear
transformation $Y=Xv(T)|\mathcal{D}$, the restriction of $Xv(T)$ to
$\mathcal{D}$, commutes with $\mathcal{W}_T$ whence $Y$ is closable.
Indeed, for any $A,B\in\mathcal{W}_T$ and $h\in\mathcal{D}(X)$ we
have $ABh\in\mathcal{D}$ and
\begin{align*}
v(T)(YA)Bh&=v(T)X(v(T)AB)h=X(v(T)^2AB)h \\
&=X(v(T)A)(v(T)B)h=(v(T)A)X(v(T)B)h \\
&=v(T)(AY)Bh,
\end{align*} where we use the facts that $v(T)AB\in H^\infty(T)$ and $X$ commutes with
$H^\infty(T)$. By the injectivity of $v(T)$, we infer from the above
equalities that $(YA)Bh=(AY)Bh$. Finally, it will be shown that if
$\{h_n\}$ is an arbitrary sequence in $\mathcal{D}(X)$ so that
$h_n\oplus(Xh_n)\rightarrow0\oplus k$ as $n\rightarrow\infty$ then
$k=0$. Indeed, the fact that $Y$ is closable shows that
$Yh_n=Xv(T)h_n=v(T)Xh_n\rightarrow v(T)k=0$ as $n\rightarrow\infty$,
and consequently $k=0$ by the injectivity of $v(T)$ again. This
finishes the proof that $X$ is closable.
\end{pf} \qed

For any operator $T$ of class $C_0$ and
$\mathcal{M}\in\mathrm{Lat}(T)$, it was shown in [3, Corollary
7.1.17] that $T$ has the property $(P)$ if and only if
$T|\mathcal{M}$ and $P_{\mathcal{M}^\perp}T|\mathcal{M}^\perp$ have
the property $(P)$. Hence we have the following theorem which
generalizes Lemma 5.6.

\begin{thm} Suppose that $T$ is an operator of class $C_0$ and
$\mathcal{M}\in\mathrm{Lat}(T)$. If $T_1=T|\mathcal{M}$ and
$T_2=P_{\mathcal{M}^\perp}T|\mathcal{M}^\perp$ then the algebra
$H^\infty(T)$ has the closability if and only if the algebras
$H^\infty(T_1)$ and $H^\infty(T_2)$ have the closability property.
\end{thm}

The proof of Lemma 5.6 combined with the conclusions in the
preceding theorem and Proposition 4.3 yield the following corollary
whose proof we omit.

\begin{cor} Suppose that $T_1$ and $T_2$ are operators of class $C_0$ so that
the algebras $H^\infty(T_1)$ and $H^\infty(T_2)$ have the
closability property. Then any densely defined linear transformation
$X$ intertwining $u(T_2)$ and $u(T_1)$ for any $u\in H^\infty$ is
closable. Moreover, $\overline{X}=B^{-1}A$ for some operator
$A\in\mathcal{I}(T_1,T_2)$ and quasiaffinity $B\in\{T_1\}'$.
\end{cor}

The preceding corollary shows that for any
$A\in\mathcal{I}(T_1,T_2)$ and injection $B\in\{T_1\}'$, $AB^{-1}$
is a closable, densely defined linear transformation intertwining
operators $u(T_2)$ and $u(T_1)$ for any $u\in H^\infty$ provided
that $H^\infty(T_1\oplus T_2)$ has the closability property. Note
that $H^\infty(T_2^*\oplus T_1^*)$ also has the closability
property. Consequently, $C^{-1}A$ is closed and densely defined for
any quasiaffinity $C\in\{T_2\}'$ since $C^{-1}A=(A^*C^{*-1})^*$ and
$A^*C^{*-1}$ is closable.

\end{document}